\documentclass{amsart}
\RequirePackage{amsthm}
\usepackage{graphicx}

\title{Some Ideals with Large Projective Dimension}
\author{Giulio Caviglia}
\address{Department of Mathematics, University of California, Berkeley, CA}
\email{caviglia@math.berkeley.edu}
\author{Manoj Kummini}
\address{Department of Mathematics, University of Kansas, Lawrence, KS}
\email{kummini@math.ku.edu}

\DeclareMathOperator{\depth}{depth}
\DeclareMathOperator{\projdim}{pd}
\DeclareMathOperator{\charact}{char}
\newcommand{\defeq}{:=}

\newtheorem{thm}{Theorem}
\newtheorem{question}[thm]{Question}
\newtheorem{cor}[thm]{Corollary}
\newtheorem{propn}[thm]{Proposition}

\begin{document}

\begin{abstract}
For an ideal $I$ in a polynomial ring over a field, a monomial support of
$I$ is the set of monomials that appear as terms in a set of minimal
generators of $I$. Craig Huneke asked whether the size of a monomial
support is a bound for the projective dimension of the ideal. We construct
an example to show that, if the number of variables and the degrees of the
generators are unspecified, the projective dimension of $I$ grows at least
exponentially with the size of a monomial support. The ideal we construct
is generated by monomials and binomials.
\end{abstract}

\maketitle

\section{Introduction}\label{sec:intro}
Let $R$ be a polynomial ring over a field $k$ and let $I \subseteq R$ be a
homogeneous ideal. Two measures of the complexity of $I$ are its projective
dimension, and its (Castelnuovo-Mumford) regularity; see
~\cite{eiscommalg}.

There have been attempts to obtain uniform bounds for projective dimension
and regularity based on numerical invariants of the ideal. Bounds on
regularity are discussed in~\cite{EisSyz05}. M. Stillman asked if there is
a bound for the projective dimension of an ideal having minimal generators
in degrees $d_1 \leq d_2 \leq \cdots \leq d_r$, when the number of
variables in the ring is not fixed. Only partial answers to this question
are known; see~\cite{Engheta05thesis}.

Related to Stillman's question, C. Huneke asked the following: is
the size of a monomial support of an ideal a bound for its
projective dimension? Here, by \emph{a monomial support} of $I$, we
mean the collection of monomials that appear as terms in a set of minimal
generators of $I$. Note that an ideal can have different monomial supports.
If $I$ is a monomial ideal, generated by $N$ monomials, then $\projdim R/I
\leq N$; this follows from the Taylor resolution of $R/I$ which has length
at most $N$~\cite[Ex. 17.11]{eiscommalg}.

We answer Huneke's question in the negative; in Sec.~\ref{sec:main}, we
construct a binomial ideal to show that the projective dimension can grow
exponentially with the size of a monomial support.  Motivated by this
example, we wonder:
\begin{question}
\label{qn:pdub}
Suppose $I \subseteq
R$ has a monomial support of $N$ monomials, counted with
multiplicity. Then what is a good upper bound for $\projdim R/I$?
\end{question}

Let $n \geq 2, d \geq 2$ be arbitrary. The ideal we construct in the
next section has a support of $2(n-1)(d-1)+n$ monomials counted with
multiplicity and projective dimension $n^d$. Using this example, we
show that for any positive integer $N$, the maximum of the
projective dimension of an ideal $I$ with a support of $N$
monomials, counted with multiplicity, is at least $2^\frac{N}{2}$.
Therefore any answer to Question~\ref{qn:pdub} should be at least
exponential. If the number of variables of $R$ is not
fixed, as in Stillman's question, the existence of any bound
is still unknown. 

Our decision of taking the multiplicity into account while counting the
monomials in the support of $I$ is only a matter of exposition. For
example, let ${m_1,\dots,m_N}$ be $N$ distinct monomials, all of the same
degree, and let $f_1,\dots,f_r$, with $f_i= \sum_{j=i}^{N}a_{ij} m_j$ and
$a_j$'s in $K,$  be a minimal system of generators for an ideal $I$. By
doing an elimination, analogous to the
one used in computing a reduced Gr\"obner basis, we can find a system of
generators $g_1,\dots,g_r$, $I=(g_1,\dots,g_r)$, such that the initial
monomial of $g_i$ does not belong to the monomial support of $g_j$ when $j
\neq i$. In this way we get a monomial support for $I$ of at most
$\sum_{i=0}^{r-1}(N-2i)=-r^2+r(N+1)$ monomials, counted with multiplicity.
The maximum value of it, as a function of $r$, is $\lfloor
(\frac{N+1}{2})^2 \rfloor$, which occurs when $r= \lfloor (N+1)/2 \rfloor$.

In general, if we have $N$ distinct monomials in a monomial support of an
homogeneous ideal $I$, then we would have at most $\lfloor
(\frac{N+1}{2})^2 \rfloor$ of them when counted with multiplicity, this is
because the above function is quadratic and the worst possible case happens
precisely when $I$ is generated by forms having the same degree.

\section{Main Example}\label{sec:main}

\begin{figure}[h]
\includegraphics{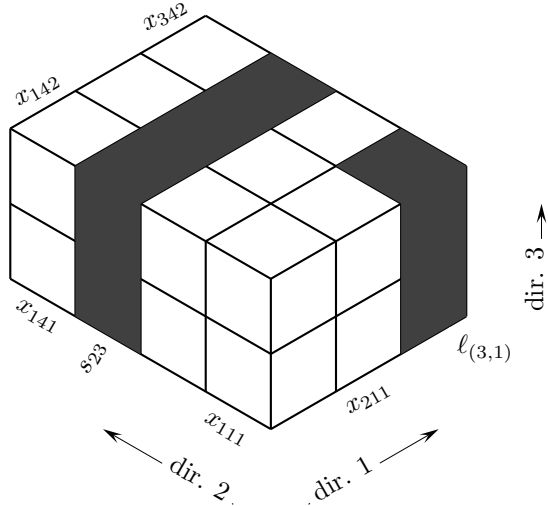}
\caption{Slices of a $3 \times 4 \times 2$ array}
\label{fig:3piped}
\end{figure}

The following example is a slight generalization of the ideal mentioned in
the introduction.
Let $d \geq 2$ and let $n_i \geq 2, \, 1 \leq i \leq d$ be positive integers.
Denote by $\mathcal I$ the index set $\{1, \cdots, n_1 \} \times \cdots
\times \{1, \cdots, n_d \}$.  Let $X \defeq \{x_{\nu} : {\nu} \in \mathcal
I\}$ be a $d$-dimensional array of variables and let $R = k[X]$. Let
\[
s_{ij} \defeq
\prod_{\substack{\nu \in \mathcal I \\ \nu_i= j}}  x_{\nu}, \;
1 \leq j\leq n_i, \; 1 \leq i \leq d.
\]
We will call $s_{ij}$ \emph{the $j$th slice in the $i$th direction}.
Fig.~\ref{fig:3piped} illustrates the above definitions for a $3 \times 4
\times 2$ array. ($\ell_.$ in the figure will defined later.)

Let $I = (s_{i 1} - s_{i j} : 2 \leq j \leq n_i, 1 \leq i \leq d-1) +
(s_{dj} : 1 \leq j \leq n_d)$. Then:
\begin{propn}
With notation as above, $\depth R/I = 0$.
\end{propn}

\begin{proof}
Write $\mathfrak{m}$ for the homogeneous maximal ideal of $R$ and let
\[
s \defeq \prod_{i = 1}^{d-1} \prod_{j = 2}^{n_i} s_{i
j}
\]
$s$ is the product of the variables not appearing in the first slices in
each of the directions $1, \cdots, d-1$. We claim that $s \in (I \colon
\mathfrak{m}) \setminus I$. Indeed, if $(I \colon \mathfrak{m}) \neq I$,
then $\mathfrak m$ is an associated prime of $R/I$, so $\depth R/I = 0$.

We first reduce the proof to the case when $\charact k = 0$, as follows.
Since $I$ is generated by monomials and binomials with $\pm 1$ as
coefficients, a Gr\"obner basis for $I$, and hence the ideal membership
problem $s \in (I \colon \mathfrak{m}) \setminus I$ are independent of the
characteristic of the field. See~\cite{eiscommalg} for the definition of
a Gr\"obner basis and the ideal membership problem. We assume, from now on,
that $\charact k = 0$.

Let $\nu \in \mathcal I$. Using the binomial relations in $I$, we can
write
\[
s \equiv \prod_{i = 1}^{d-1} s_{i 1} \cdots \widehat{s_{i \nu_i}} \cdots
s_{i n_i} \mod I
\]
where $\widehat \cdot$ denotes omitting the variable from the
product.
Consider the slice $s_{d\nu_d} =
\prod_{\substack{\mu \in \mathcal I \\ \mu_d= \nu_d}}  x_{\mu}$. If $\mu
\neq \nu \in \mathcal I$ is such that $\mu_d = \nu_d$, then there exists $1
\leq i \leq d-1$ such that $\mu_i \neq \nu_i \implies x_\mu \vert (s_{i 1}
\cdots \widehat{s_{i \nu_i}} \cdots s_{i n_i}) \implies s_{d\nu_d} \vert
((\prod_{i=1}^{d-1} s_{i 1} \cdots \widehat{s_{i \nu_i}} \cdots s_{i
n_i})x_\nu) \implies ((\prod_{i=1}^{d-1} s_{i 1} \cdots \widehat{s_{i
\nu_i}} \cdots s_{i n_i})x_\nu) \in I$. Hence $s \in (I \colon
\mathfrak{m})$.

Let $A$ be the tableau
\[
\begin{array}{lll}
a_{11} & \cdots&  a_{1n_1}\\
a_{21} & \cdots&  a_{2n_2}\\
&\ddots&\\
a_{(d-1)1} & \cdots&  a_{(d-1)n_{(d-1)}}
\end{array}
\]
of non-negative integers. We use \emph{tableau} loosely here; we only mean
that the rows of $A$ possibly have different number of elements. For such
a tableau $A$,  we say it \emph{satisfies row condition $(c_1, \cdots,
c_t)$} if the sum of the elements on the $i$th row is $c_i-1$.

Let $\mathcal P \defeq \{1, \cdots, n_1 \} \times \cdots \times \{1,
\cdots, n_{d-1}\}$. For each $p \in \mathcal P$, we define a monomial
\[
\ell_p \defeq
\prod_{\substack{\nu \in \mathcal I \\ \nu_i= p_i, 1 \leq i \leq d-1}}
x_{\nu}, \]
See Fig.~\ref{fig:3piped} for an illustration of $\ell_{(3,1)}$ in the $3
\times 4 \times 2$ case. Further, write $\vert p \vert_A$ for
$\sum_{i=1}^{d-1}{a_{ip_i}}$

Let
\begin{equation}\label{eqn:masterpoly}
F = \sum_{A : A \text{ satisfies }(n_1, \cdots, n_{d-1})} \left\{\prod_{p \in
\mathcal P} \frac{1}{(\vert p \vert_A!)^{n_d}} \ell_p^{\vert p \vert_A}
\right\}
\end{equation}
We let $R$ act on itself by partial differentiation with respect to the
variables. We show below that, under this action, $s \not \in (0 \colon_R
F)$ while $I \subseteq (0 \colon_R F)$ from which we conclude that $s \not
\in I$, thus proving the proposition.

For any tableau $A$ that satisfies the row condition $(n_1, \cdots,
n_{d-1})$, write $\tau_A$ for the corresponding monomial term that appears
in $F$ (see \eqref{eqn:masterpoly}).  Let
\[
A_s \defeq
\begin{array}{llll}
0 & 1 & \cdots& 1\\
0 & 1 & \cdots& 1\\
&&\ddots&\\
0 & 1 & \cdots&  1\\
\end{array}
\]
Then $s = \alpha \tau_{A_s}$ for some non-zero rational number
$\alpha$ If $A \neq A_s$, then $s$ contains a
variable that $\tau_A$ does not contain, so $s \circ F = s \circ \tau_{A_s}
= 1$. Hence $s \not \in (0 \colon_R F)$.

For any $1 \leq j \leq n_d$, $s_{dj} \circ F = 0$. For, any $A$ that
appears in the summation of \eqref{eqn:masterpoly} has at least one $p_A
\in \mathcal P$ such that $\vert p_A \vert_A = 0$. Hence the variables in
$\ell_{p_A}$ do not appear in $\tau_A$. However, $s_{dj}$ contains
one such variable, and hence, $s_{dj} \circ \tau_A = 0 \implies s_{dj}
\circ F = 0$.

Observe that any slice $s_{ij}, 1 \leq i \leq d-1, 1 \leq j \leq n_i$ can
be written as a product of $\ell_p, p \in \mathcal P$ as follows:
\[
s_{ij} =
\prod_{\substack{\nu \in \mathcal I \\ \nu_i= j}}  x_{\nu} =
\prod_{\substack{1 \leq \nu_{i'} \leq n_{i'} \\ 1 \leq i' \leq d-1 \\
\nu_i= j}}
\underbrace{
\left[ \prod_{\substack{1 \leq \nu_d \leq n_d}} x_{\nu}
\right]
}_{\ell_{(\nu_1, \cdots, \nu_{d-1})}} =
\prod_{\substack{p \in \mathcal P \\ p_i = j}} {\ell_p} \]
Let $\mathcal P_{ij} = \{ p \in \mathcal P : p_i = j\}$. Then $s_{ij}
\circ F = \left(\prod_{p \in \mathcal P_{ij}} \ell_p \right) \circ F$.
Therefore to differentiate with respect to $s_{ij}$, we may differentiate
with respect to all $\ell_p, p \in \mathcal P$, sequentially.

Let $1 \leq i \leq d-1, 1 \leq j \leq n_i$ and $q \in \mathcal P_{ij}$. Then
\begin{align*}
\ell_q \circ F & =
\ell_q \circ \sum_{A : A \text{ satisfies }(n_1, \cdots, n_{d-1})}
\left\{\prod_{p \in \mathcal P} \frac{1}{(\vert p \vert_A!)^{n_d}}
\ell_p^{\vert p \vert_A} \right\}  \\
& =
\sum_{A : A \text{ satisfies }(n_1, \cdots, n_{d-1})} \left\{
\frac{(\vert q \vert_A)^{n_d}}{(\vert q \vert_A!)^{n_d}} \ell_q^{(\vert q
\vert_A - 1)}
\prod_{\substack{p \in \mathcal P\\ p \neq q}}
\frac{1}{(\vert p \vert_A!)^{n_d}} \ell_p^{\vert p \vert_A}
\right\}\\
\end{align*}
Therefore,
\begin{equation}\label{eqn:partialwrtsij}
s_{ij} \circ F =
\sum_{A : A \text{ satisfies }(n_1, \cdots, n_{d-1})} \left\{
\prod_{p \in \mathcal P_{ij}}
\frac{(\vert p \vert_A)^{n_d}}{(\vert p \vert_A!)^{n_d}} \ell_p^{(\vert p
\vert_A - 1)}
\prod_{p \not \in \mathcal P_{ij}}
\frac{1}{(\vert p \vert_A!)^{n_d}} \ell_p^{\vert p \vert_A}
\right\}\\
\end{equation}
We can write $\{A : A \text{ satisfies }(n_1, \cdots, n_{d-1}) \} = \{ A :
a_{ij} = 0 \} \bigcup \{ A : a_{ij} \neq 0 \}$. Every row of $A$ contains
at least one zero. If $a_{ij} = 0$, then there is a $p \in \mathcal P_{ij}$
such that $\vert p \vert_A = 0$. Therefore there is no contribution from
those $A$ with $a_{ij} = 0$ in the RHS of \eqref{eqn:partialwrtsij}.
Moreover, $a_{ij} \neq 0 \implies \vert p \vert_A \neq 0$. Hence
\[
s_{ij} \circ F =
\sum_{\substack{A \text{ satisfies }(n_1, \cdots, n_{d-1}) \\ a_{ij} \neq
0}} \left\{ \prod_{p \in \mathcal P_{ij}}
\frac{1}{[(\vert p \vert_A - 1)!]^{n_d}} \ell_p^{(\vert p
\vert_A - 1)}
\prod_{p \not \in \mathcal P_{ij}}
\frac{1}{(\vert p \vert_A!)^{n_d}} \ell_p^{\vert p \vert_A}
\right\}
\]
There is a 1-1 correspondence between
$\{A : A \text{ satisfies }(n_1, \cdots, n_{d-1}), a_{ij} \neq
0\}$ and  $\{A : A \text{ satisfies }(n_1, \cdots, n_i - 1, \cdots,
n_{d-1}) \}$. Using this we can write
\begin{equation} \label{eqn:partialwrtsijnew}
s_{ij} \circ F =
\sum_{A \text{ satisfies }(n_1, \cdots, n_i - 1, \cdots, n_{d-1})}
\left\{\prod_{p \in \mathcal P} \frac{1}{(\vert p \vert_A!)^{n_d}}
\ell_p^{\vert p \vert_A} \right\}
\end{equation}
Note that this representation of $s_{ij} \circ F$ is independent of $j$;
hence $(s_{i1}-s_{ij})\circ F =  0$ for all $1 \leq i \leq d-1$ and $2
\leq j \leq n_i$. Hence $I \subseteq (0 \colon_R F)$.
\end{proof}

It now follows from the Auslander-Buchsbaum formula (see,
\textit{e.g},~\cite[Theorem 19.9]{eiscommalg} that
\begin{cor}
With notation as above, $\projdim R/I = n_1 \cdots n_d$.
\qed
\end{cor}

Parenthetically, we note that the ideal we construct has $n_i-1$ generators
of degree $n_1 \cdots \widehat{n_i} \cdots n_d$, for $1 \leq i \leq d-1$
and $n_d$ generators of degree $n_1 \cdots n_{d-1}$.

Consider the case when $n_1 = \cdots = n_d = n$. Then the ideal is generated
by $(n-1)(d-1)$ binomials and $n$ monomials, and, hence, has a monomial
support of $2(n-1)(d-1)+n$.

\begin{cor}
Any upper bound for projective dimension of an ideal supported on $N$
monomials counted with multiplicity is at least $2^{N/2}$.
\end{cor}

\begin{proof}
Given a positive integer $N$, choose $n=2$ variables in each of $d =
\frac{N}{2}$ dimensions, and construct $R$ and $I$ as above. Then $\projdim
R/I = 2^{N/2}$.
\end{proof}

\section*{Acknowledgements}
We thank Profs D. Eisenbud, C. Huneke and B. Sturmfels for helpful
discussion.

\providecommand{\bysame}{\leavevmode\hbox to3em{\hrulefill}\thinspace}
\providecommand{\MR}{\relax\ifhmode\unskip\space\fi MR }
% \MRhref is called by the amsart/book/proc definition of \MR.
\providecommand{\MRhref}[2]{%
  \href{http://www.ams.org/mathscinet-getitem?mr=#1}{#2}
}
\providecommand{\href}[2]{#2}

\end{document}